\documentclass[11pt]{article}
\usepackage{amssymb}

\newcommand{\proof}{\noindent{\bf Proof. }}
\newcommand{\fg}{{\mathfrak g}}
\newcommand{\fG}{{\mathfrak G}}

\newcommand{\fk}{{\mathfrak k}}
\newcommand{\fz}{{\mathfrak z}}
\newcommand{\fh}{{\mathfrak h}}
\newcommand{\fs}{{\mathfrak s}}
\newcommand{\fI}{{\mathfrak I}}

\newcommand{\ad}{\,\mbox{\rm ad}}

\newcommand{\qed}{\hfill $\square$}

\newcommand{\C}{\mathbb{C}}
\newcommand{\R}{\mathbb{R}}

\newtheorem{theor}{Theorem}[section]
\newtheorem{prop}[theor]{Proposition}
\newtheorem{coro}[theor]{Corollary}
\newtheorem{lemma}[theor]{Lemma}
\newtheorem{rema}{Remark}[section]
\newtheorem{defi}{Definition}[section]
\newtheorem{example}{Example}[section]

\begin{document}

\title{On Pseudo-Hermitian Quadratic Nilpotent Lie Algebras}

\author{Mustapha Bachaou, Ignacio {Bajo\footnote{Corresponding author}$\,$} and Mohamed Louzari} 

\maketitle

\begin{abstract} We study nilpotent Lie algebras endowed with a complex structure and a  quadratic structure which is pseudo-Hermitian for the given complex structure. We propose several methods to construct such Lie algebras and describe a method of  double extension by  planes to get an inductive description of all of them. As an application, we give a complete classification of nilpotent quadratic Lie algebras where the metric is Lorentz-Hermitian  and  we fully classify all nilpotent pseudo-Hermitian quadratic Lie algebras up to dimension 8 and their inequivalent pseudo-Hermitian metrics.
\end{abstract}



\section*{Introduction}

A pseudo-Hermitian quadratic Lie algebra is a triple $(\fg,J,\varphi)$ where $\fg$ is a Lie algebra, $J$ is a complex structure on $\fg$ and $\varphi$ is a non-degenerate ad-invariant bilinear form on $\fg$ which is compatible with $J$ in the sense that $\varphi(Jx,Jy)=\varphi(x,y)$ for all $x,y\in\fg$. Geometrically, $\fg$ is the Lie algebra of a Lie group $G$ endowed with a left-invariant complex structure and a bi-invariant pseudo-Hermitian metric. Further, it was seen in \cite{Andrada} that the complex structure on such Lie algebras defines a classical $r$-matrix so that they give rise to  Lie bialgebras of complex type and, as a consequence, to certain Poisson Lie groups whose duals are complex Lie groups.  

Even though the complex structure on a pseudo-Hermitian quadratic Lie algebra cannot be neither bi-invariant nor abelian, one may see that the class of such Lie algebras is a large one and includes nilpotent, solvable and non-solvable examples; one may even find examples of pseudo-Hermitian quadratic Lie algebras that are perfect. In this paper we will focus our attention in the case of nilpotent ones. 

Medina and Revoy \cite{MR} gave an inductive method to construct all solvable quadratic Lie algebras. The method is commonly known as the quadratic double extension by a line method and consists, essentially, in describing each solvable quadratic Lie algebra of dimension $n$ as the Lie algebra obtained from a $n-2$-dimensional quadratic Lie algebra after a central extension by a 1-dimensional ideal and a suitable semidirect product by another 1-dimensional Lie algebra. When working with Lie algebras provided with complex structures one wants the Lie algebra of dimension $n-2$ to be also endowed with a complex structure and, in most of the cases, this situation cannot be described by a single double extension by a line. Thus, in order to guarantee the existence of a complex structure in the Lie algebra of lower dimension, one needs to consider a couple of successive double extensions. In the nilpotent case this can be formulated as a certain generalized double extension in which the central extension is done by a 2-dimensional $J$-invariant subspace and one applies two semidirect products. As a rule, the two derivations with which one does the semidirect products do not span a Lie subalgebra of the large pseudo-Hermitian quadratic Lie algebra since they may have a non-zero bracket lying in the small Lie algebra. In our Proposition \ref{pR1} we give a detailed description of this modified double extension method and in Theorem \ref{mainth} we show that every indecomposable nilpotent pseudo-Hermitian quadratic Lie algebra can be constructed in that way. A simple application of such a method let us describe all Lie algebras that  can be endowed with Lorentz-Hermitian quadratic metrics. 

Using our double extension method,  sometimes called pHQ-double extension within the text, we give a complete classification of nilpotent pseudo-Hermitian quadratic Lie algebras. It should be remarked that, although quadratic nilpotent Lie algebras up to dimension 10 which are indecomposable (as quadratic Lie algebras)  have been classified in \cite{IK}, one cannot easily deduce our classification from the results therein. On the one hand, indecomposability as a quadratic Lie algebra is not the same as indecomposability as a pseudo-Hermitian quadratic Lie algebra for one may find examples of pseudo-Hermitian quadratic Lie algebras that cannot be decomposed as a sum of two pseudo-Hermitian quadratic Lie algebras but can be split as the sum of two quadratic ones; this is the case, for instance, of the Lie algebra ${\mathcal W}_3\oplus \R$, where ${\mathcal W}_3$ is as defined in \cite{FS}, which is obviously decomposable as quadratic Lie algebra but indecomposable as pseudo-Hermitian quadratic Lie algebra. On the second hand, a quadratic Lie algebra might admit non-equivalent complex structures or pseudo-Hermitian quadratic metrics. In Theorem \ref{c8} we give the complete list of inequivalent  pseudo-Hermitian quadratic structures on Lie algebras up to dimension 8. Precisely, we show that there are only six inequivalent indecomposable pseudo-Hermitian quadratic Lie algebras up to dimension 8.

\section{Definitions, examples and first results}

From now on, unless otherwise stated, all the Lie algebras considered in the paper are real Lie algebras.

\begin{defi}{\em 
A {\it complex structure} on an  Lie algebra $\mathfrak{g}$ is a map $J \in \mathfrak{gl(g)}$ such that
$J^2 = -I$, where $I$ stands for the identity map, and the equality
\begin{equation}\label{nijenhuis}
N_{J}(x,y)=[x, y] + J[Jx, y] + J[x, Jy] - [Jx, Jy] = 0 
\end{equation}
holds for all $x, y \in \mathfrak{g}$. 

The complex structure $J$ is said to be {\it abelian} if $[Jx,Jy]=[x,y]$ for every $x,y\in\fg$ whereas it is said to be {\it bi-invariant} when 
 $[Jx,y]=J[x,y]$ holds for every $x,y\in\fg$.

Two Lie algebras endowed with complex structures $(\fg_1,J_1)$, $(\fg_2,J_2)$ are said to be {\it holomorphically isomorphic} if there exists an isomorphism of Lie algebras $\psi:\fg_1\to\fg_2$ such that $J_2\circ\psi=\psi\circ J_1$.
}
\end{defi}

\begin{rema}\label{rem1}{\em Let us point out a couple of interesting facts.
\begin{enumerate}
\item[(1)] The so-called Nijenhuis tensor $N_{J}(x,y)=[x, y] + J[Jx, y] + J[x, Jy] - [Jx, Jy]$  is clearly antisymmetric and verifies the following identities:
\begin{equation}\label{propnijen}
N_J(Jx,Jy)=-N_J(x,y),\quad N_J(Jx,y)=N_J(x,Jy)=-JN_J(x,y),
\end{equation}
for all $x,y\in\fg$.
\item[(2)] Denote by $\fg^{1,0}=\mbox{Ker}(J-iI)$ and $\fg^{0,1}=\mbox{Ker}(J+iI)$ the corresponding eigenspaces of eigenvalues $i$ and $-i$ in the complexification $\fg^\C$ of $\fg$, this is to say:
$$
 \fg^{1,0}=\{x-iJx\ :\ x\in\fg\},\quad
\fg^{0,1}=\{x+iJx\ :\ x\in\fg\}.$$
One has that equality (\ref{nijenhuis}) is equivalent to $\fg^{1,0}$ (or $\fg^{0,1}$) being a complex Lie subalgebra of $\fg^\C$. It is well known that $J$ is abelian if and only if $\fg^{1,0}$ is an abelian Lie algebra and that $J$ is bi-invariant if and only if $(\fg, J)$ and $(\fg^{1,0}, i)$ are holomorphically isomorphic (and, hence, $\fg$ can be considered as a complex Lie algebra).

Moreover, for an arbitrary complex structure $J$ on a Lie algebra $\fg$, 
denote by $\fg_J$ the vector space $\fg$ endowed with a new bracket defined by
$$[x,y]_J=[Jx,y]+[x,Jy],\quad\ x,y\in\fg.$$
It can be seen that $[\cdot,\cdot]_J$ is actually a Lie bracket and that there is a natural holomorphic isomorphism $\psi:\fg_J\to \fg^{1,0}$ defined by $\psi (x)=Jx+ix.$

It should be remarked that the map $\eta:\fg_J\to \fg^{0,1}$ defined by $\eta(x)=Jx-ix$ is also an isomorphism of real Lie algebras. However, $\eta$ is not holomorphic but anti-holomorphic since
$$\eta(Jx)=-x-iJx=-i(Jx-ix)=-i\eta(x),\quad\ x\in\fg.$$
\item[(3)] Notice that every Lie algebra with a complex structure $(\fg,J)$ has the structure of complex vector space defined by the multiplication by scalars given by:
$$(\alpha+i\beta)\cdot x=\alpha x+\beta Jx,\quad \alpha,\beta\in\R,\ x\in\fg.$$
As a consequence, every Lie algebra admitting a complex structure must have even dimension. Actually, this is also true for every Lie algebra which admits a linear endomorphism $J$ such that $J^2=-I$ (even if the Nijenjuis tensor is not identically zero).
\end{enumerate}
}
\end{rema}

\begin{defi} {\em Let $(\fg, J)$ be a Lie algebra with a complex structure.
If $\mathfrak{g}$ is also endowed with a non-degenerate symmetric bilinear map 
$
\varphi : \mathfrak{g}\times \mathfrak{g} \to\mathbb{R} $
  verifying the condition $\varphi ( Jx, Jy) = \varphi ( x, y),$ for all $x, y \in \mathfrak{g}$, then $(\mathfrak{g}, J, \varphi )$ is said to be a {\it pseudo-Hermitian Lie algebra} and the pair $(J, \varphi)$ a {\it pseudo-Hermitian structure} on $\mathfrak{g}$.
	
	For a pseudo-Hermitian Lie algebra $(\mathfrak{g}, J, \varphi )$ one defines the {\it fundamental 2-form} or {\it K\"ahler form} by
	$$\omega(x,y)=\varphi(x,Jy),\quad\ x,y\in\fg.$$
	
	\smallskip

	It is an easy exercise to see that the  map $(\cdot |\cdot):\fg\times\fg\to\C$ defined by
$$(x|y)=\varphi( x,y)+i\omega(x,y),\quad\ x,y\in\fg$$
is actually a Hermitian form on the complex vector space $\fg$ defined as in Remark \ref{rem1} (3), this is to say that $(\cdot |\cdot)$ is a sesquilinar form  such that $(y|x)=\overline{(x|y)}$ for every $x,y\in\fg$.}\end{defi}

\begin{rema}{\em It should be remarked that the signature of a pseudo-Hermitian metric $\varphi$ on a $2n$-dimensional Lie algebra is always of the form $\mbox{sig}(\varphi)=(2r,2s)$. We will use the term {\it Hermitian metric} exclusively for the definite positive case $\mbox{sig}(\varphi)=(0,2n)$ and we will call a {\it Lorentz-Hermitian metric} the case of index 2, this is to say, $\mbox{sig}(\varphi)=(2,2n-2)$.

Also, notice that, since $J^2=-I$, the condition $\varphi ( Jx, Jy) = \varphi ( x, y),$ is equivalent to $\varphi ( Jx, y) = -\varphi ( x, Jy),$ for all $x,y\in\fg$.}
\end{rema}

\begin{defi}{\em Let $(\mathfrak{g}_1, J_1, \varphi_1 )$, $(\mathfrak{g}_2, J_2, \varphi_2 )$ be two pseudo-Hermitian Lie algebras. We will say that they are {\it equivalent} if there exists a holomorphic isomorphism $\psi :\fg_1\to\fg_2$ which is also an isometry, this is to say,  $\varphi_2(\psi x,\psi y)=\varphi_1(x,y)$ holds for all $x,y\in\fg_1.$}
\end{defi}

\begin{defi}{\em Let $\fg$ be a Lie algebra and $\varphi:\fg\times\fg\to\R$ a non-degenerate symmetric bilinear form. If for every $x,y,z\in\fg$ it holds that
$$\varphi([x,y],z)+\varphi(y,[x,z])=0,$$
then we say that $\varphi$ is a {\it quadratic structure} or a {\it quadratic metric} on $\fg$. In such case, we also say that $(\fg ,\varphi)$ is a  {\it quadratic Lie algebra}. In many references, quadratic Lie algebras are also called {\it metric Lie algebras} \cite{IK}.

If a quadratic Lie algebra $(\fg,\varphi)$ also admits a complex structure $J$ such that $(\mathfrak{g}, J, \varphi )$ is pseudo-Hermitian, we will say that $(\mathfrak{g}, J, \varphi )$ is a {\it pseudo-Hermitian quadratic Lie algebra.}}
\end{defi}

\medskip

Our next result shows that  abelian complex structures do not exist in quadratic Lie algebras unless $\fg$ itself is abelian (even if the metric is not pseudo-Hermitian).

\medskip

\begin{prop} If $(\mathfrak{g},\varphi )$ is a non-abelian quadratic Lie algebra, then $\fg$ cannot admit an abelian complex structure.
\end{prop}
\proof Suppose that $(\mathfrak{g},\varphi )$ is a quadratic Lie algebra and $J$ an abelian complex strtucture on $\fg$. Then we have
\begin{eqnarray*} 
 \varphi( [Jx,Jy],z) &=& -\varphi( Jy,[Jx,z]) =  
 \varphi(Jy,[x,Jz] )  =-\varphi( Jy,[Jz,x] )=\varphi( [Jz,Jy],x )\\&=&\varphi( [z,y],x )=-\varphi( [y,z],x )=\varphi( z,[y,x ])=\varphi( z,[Jy,Jx ])\\
&=&-\varphi( z,[Jx,Jy ])=
- \varphi([Jx,Jy], z ), 
\end{eqnarray*}
	for all $x,y,z\in\fg$, which clearly shows that $\fg$ must be abelian because $\varphi$ is non-degenerate.\qed 
	
\begin{rema}{\em In contrast with the proposition above, non-abelian quadratic Lie algebras may admit bi-invariant complex structures. For instance, if $\fg=\fs_\R$ is the underlying real Lie algebra of a complex semisimple Lie algebra $\fs$ and $\varphi$ denotes the real part of the Killing form ${\mathcal K}$ on $\fs$, it is easy to prove that $(\fg,\varphi)$ is quadratic and that the map $J:\fg\to\fg$ defined by $Jx=ix$ for all $x\in\fs$ is a bi-invariant complex structure on $\fg$. Notice, however, that the triple $(\fg, J,\varphi)$ is not pseudo-Hermitian since
$$\varphi(Jx,Jy)=\mbox{Re}({\mathcal K}(ix ,iy))=- \mbox{Re}({\mathcal K}(x ,y))=-\varphi(x,y),$$
for every $x,y\in\fs$. These kind of metrics verifying $\varphi(Jx,Jy)=-\varphi(x,y)$ are usually called {\it Norden} metrics (see \cite{bajo1} and references therein).

In fact, the following result is well known (see, for instance, \cite[Prop. 9]{Andrada}):}
\end{rema}

\begin{prop} Let $(\fg, J,\varphi)$ be a pseudo-Hermitian quadratic Lie algebra. The complex structure $J$ is bi-invariant if and only if $\fg$ is an abelian Lie algebra.
\end{prop}

\begin{rema}{\em Although there are no non-abelian examples of pseudo-Hermitian quadratic Lie algebras with abelian or bi-invariant complex structures, the class of such Lie algebras is a wide one.

In \cite[Th. 24]{Andrada} it is seen that,  up to dimension 6, every quadratic Lie algebra $(\fg,\varphi)$ with $\mbox{sig}(\varphi)=(2r,2s)$ admits a $\varphi$-skewsymmetric complex structure $J$ and, therefore, $(\fg, J,\varphi)$ is pseudo-Hermitian quadratic.

We will next give some constructions of pseudo-Hermitian quadratic Lie algebras.}
\end{rema}

\medskip

Our first construction uses the concept of T$^*$-extension defined by Bordemann in \cite{Borde}. 

\medskip

\begin{prop}\label{T*} Let $(\fg,J)$ be a Lie algebra with a complex structure and suppose that there exists a  2-cocycle $\theta\in Z^2(\fg,\fg^*)$ for the Chevalley-Eilenberg complex defined by the coadjoint representation of $\fg$ such that for all $x,y,x\in\fg$ the  following equality  holds:
\begin{equation}\label{compJ} \theta(x,y)z=\theta(Jx,Jy)z+\theta(Jy,Jz)x+\theta(Jz,Jx)y .\end{equation}

Define on the vector space $T^*_\theta\fg=\fg\oplus\fg^*$ the bilinear form $\varphi (x+f,y+g)=f(y)+g(x)$  and the bracket 
$$[x+f,y+g]_T=[x,y]+\theta(x,y)+f\circ\ad(y)-g\circ\ad(x),$$
for all $x,y\in\fg$, $f,g\in\fg^*$, where $\ad (x)y=[x,y]$, and let $J_T$ be defined by $J_T(x+f)=Jx-f\circ J.$

The triple $(T^*_\theta\fg,J_T,\varphi)$ is a pseudo-Hermitian quadratic Lie algebra.
\end{prop}
\proof First, notice that for every $x,y,z\in\fg$ we have
$$\theta(y,z)x=\theta(Jy,Jz)x+\theta(Jz,Jx)y+\theta(Jx,Jy)z=\theta(x,y)z,$$
which proves the cyclic condition of \cite[Lemma 3.1]{Borde}. This implies that
 $(T^*_\theta\fg,\varphi)$ is a quadratic Lie algebra by the results in \cite{Borde}. So, we only have to prove that $J_T$ is a complex structure on $T^*_\theta\fg$ and that it is skewsymmetric with respect to $\varphi$. 

This last fact is clear since
\begin{eqnarray*}\varphi(J_T(x+f),y+g)&=&\varphi(Jx-f\circ J),y+g)=g(Jx)-f(Jy)=
-\varphi(x+f,Jy-g\circ J)\\&=&-\varphi(x+f,J_T(y+g)),\end{eqnarray*}
for every $x,y\in\fg$, $f,g\in\fg^*$. 

To see that $J_T$ is a complex structure, first notice that $J_T^2(x+f)=J^2x+f\circ J^2=-x-f$, for all $x\in\fg,f\in\fg^*$.
Let us see that $N_{J_T}$ is identically zero. For $f,g\in\fg^*$ we clearly have
$$
N_{J_T}(f,g)=[f,g]_T-J_T[f\circ J,g]_T-J_T[f,g\circ J]_T-[f\circ J,g\circ J]_T=0.$$
Now, choosing $x\in\fg,g\in\fg^*$ we get
\begin{eqnarray*}N_{J_T}(x,g)&=&[x,g]_T+J_T[Jx,g]_T-J_T[x,g\circ J]_T+[Jx,g\circ J]_T\\&=&
-g\circ\ad(x)+g\circ\ad(Jx)\circ J-g\circ J\circ\ad(x)\circ J-g\circ J\circ \ad(Jx)\in\fg^*,\end{eqnarray*}
but then $N_{J_T}(x,g)=0$ because
$$N_{J_T}(x,g)(y)=-g([x,y]-[Jx,Jy]+J[x,Jy]+J[Jx,y])=-g(N_J(x,y))=0.$$
Finally, for $x,y\in\fg$ one gets
\begin{eqnarray*}
N_{J_T}(x,y)&=&[x,y]_T+J_T[Jx,y]_T+J_T[x, Jy]_T-[Jx,Jy]_T\\
&=&N_J(x,y)+\theta(x,y)-\theta(Jx,y)\circ J-\theta(x,Jy)\circ J-\theta(Jx,Jy).\end{eqnarray*}
But, using the hypothesis on $\theta$ and its cyclicity, we have for all $z\in\fg$ that 
\begin{eqnarray*}0&=&\theta(x,y)z-\theta(Jz,Jx) y-\theta(Jy, Jz)x-\theta(Jx,Jy)z\\&=&\theta(x,y)z-\theta(Jx,y) Jz-\theta(x,Jy) Jz-\theta(Jx,Jy)z,\end{eqnarray*}
so that
$N_{J_T}(x,y)=N_J(x,y)=0,$
which completes the proof.\qed

\begin{example}{\em Using Proposition \ref{T*} we can construct a lot of non-trivial examples of pseudo-Hermitian quadratic Lie algebras. Let us mention some interesting cases: 
\begin{enumerate}
\item[(1)] If one considers $\theta=0$ then the Lie algebra $T^*_0\fg$ is actually the cotangent algebra $T^*\fg$ \cite[Example 11]{Andrada}. Notice that this method let us construct a pseudo-Hermitian quadratic Lie algebra starting from any pair $(\fg,J)$ of a Lie algebra with a complex structure.
\item[(2)] A  interesting particular case of the construction above is when one considers $\fg$ the underlying real Lie algebra of a semisimple complex Lie algebra with the natural bi-invariant complex structure defined by $Jx=ix$, for $x\in\fg$. In this case the cotangent Lie algebra $T^*\fg$ provides an example of pseudo-Hermitian quadratic Lie algebra which is perfect (this is to say, $[T^*\fg,T^*\fg]_T=T^*\fg$) and, hence, centerless.

Obvously, the same construction can be done for a rank 2 compact semisimple Lie algebra endowed with a complex structure \cite{Samelson} to get that its cotangent Lie algebra is a perfect pseudo-Hermitian quadratic Lie algebra. 
\item[(3)] One may find examples in which $\theta\ne 0$. For instance, let  $\fk$ be the direct sum of the 3-dimensional Heisenberg algebra and $\R$, this is to say $\fk=\R\mbox{-span}\{x_1,x_2,x_3,x_4\}$ with the only non trivial bracket $[x_1,x_2]=x_3$. Such a Lie algebra is the well-known algebra underlying the Kodaira-Thurston manifold \cite{Koda,Thurs} and it is well known that the linear map $J$ defined by $Jx_1=x_2$, $Jx_2=-x_1$, $Jx_3=x_4$, $Jx_4=-x_3$ is an abelian complex structure on $\fg$. In our next result we give a complete description of all non-zero admissible cocycles $\theta$ and of the corresponding Lie algebras $T^*_\theta\fk$. Those examples will play a relevant role in our classification of 8-dimensional pseudo-Hermitian quadratic Lie algebras given in Section \ref{classif}.\end{enumerate}}
\end{example} 

\begin{prop} \label{kodaira} Let $(\fk, J)$ be the Lie algebra $\fk=\R\mbox{\rm{-span}}\{x_1,x_2,x_3,x_4\}$  with the only non trivial bracket $[x_1,x_2]=x_3$ and the abelian complex structure given by $Jx_1=x_2$, $Jx_2=-x_1$, $Jx_3=x_4$, $Jx_4=-x_3$.

\begin{enumerate}
\item[{\rm (1)}] Every cyclic 2-cocycle $\theta\in Z^2(\fk,\fk^*)$ is a linear combination of the the cyclic cocycles $\theta_i$, $1\le i\le 4$ defined, up to skewsymmetry, by the following  non-zero values:
\begin{eqnarray*} & & \theta_1(x_1,x_2)=x_3^*,\quad \theta_1(x_1,x_3)=-x_2^*,\quad \theta_1(x_2,x_3)=x_1^*,\\
& & \theta_2(x_1,x_2)=x_4^*,\quad \theta_2(x_1,x_4)=-x_2^*,\quad \theta_2(x_2,x_4)=x_1^*,\\
& & \theta_3(x_1,x_3)=x_4^*,\quad \theta_3(x_1,x_4)=-x_3^*,\quad \theta_3(x_3,x_4)=x_1^*,\\
& & \theta_4(x_2,x_3)=x_4^*,\quad \theta_4(x_2,x_4)=-x_3^*,\quad \theta_4(x_3,x_4)=x_2^*.\end{eqnarray*}
Further, any such a 2-cocycle $\theta$ verifies the compatibility condition of eq. (\ref{compJ}).
\item[{\rm (2)}] The pseudo-Hermitian $T^*$-extension of $(\fk, J)$ by means of the cocycle $\theta=\sum_{i=1}^4\alpha_i\theta_i$ is the Lie algebra  spanned by a basis $\{x_1,x_2,x_3,x_4,x_1^*,x_2^*,x_3^*,x_4^*\}$ with  nontrivial brackets
\begin{eqnarray*}
& & [x_1,x_2]=x_3+\alpha_1x_3^*+\alpha_2x^*_4,\quad [x_1,x_3]=-\alpha_1x^*_2+\alpha_3x_4^*,\\
& & [x_1,x_4]=-\alpha_2x_2^*-\alpha_3x^*_3,\quad [x_2,x_3]=\alpha_1x^*_1+\alpha_4x_4^*,\\
& & [x_2,x_4]=\alpha_2x_1^*-\alpha_4x^*_3,\quad [x_3,x_4]=\alpha_3x^*_1+\alpha_4x_2^*,\\
& & [x_3^*,x_1]=x_2^*,\quad [x_3^*,x_2]=-x^*_1\end{eqnarray*} 
 and the metric and complex structure defined by
$$\varphi(x_i,x_j^*)=\delta_{ij},\,\, \varphi(x_i,x_j)=\varphi(x_i^*,x_j^*)=0,\quad Jx_1=x_2,\,\, Jx_3=x_4,\,\, Jx_1^*=x_2^*,\,\, Jx_3^*=x_4^*.$$  
\end{enumerate}
\end{prop}
\proof We will only prove the first statement since (2) follows at once from the construction of $T^*_\theta \fk$ as in Proposition \ref{T*}.

Let $\theta:\fk\times\fk\to\fk^*$ be an skewsymmetric map such that the cyclic condition $\theta(x,y)z=\theta(y,z)x$ is verified for all $x,y,z\in\fk$. Consider $a_i\in\R$, $1\le i\le 12$ such that 
\begin{eqnarray*}& & \theta(x_1,x_2)=a_1 x^*_3+a_2 x^*_4,\quad \theta(x_1,x_3)=a_3 x^*_2+a_4 x^*_4 \\
& & \theta(x_1,x_4)=a_5 x^*_2+a_6 x^*_3,\quad \theta(x_2,x_3)=a_7 x^*_1+a_8 x^*_4\\
& & \theta(x_2,x_4)=a_9 x^*_1+a_{10} x^*_3,\quad \theta(x_3,x_4)=a_{11} x^*_1+a_{12} x^*_2.
\end{eqnarray*}
From the cyclicity of $\theta$ we have
\begin{eqnarray*}& &a_1=\theta(x_1,x_2)x_3=\theta(x_2,x_3)x_1=a_7=\theta(x_3,x_1)x_2=-a_3\\
& &a_2=\theta(x_1,x_2)x_4=\theta(x_2,x_4)x_1=a_9=\theta(x_4,x_1)x_2=-a_5\\
& &a_4=\theta(x_1,x_3)x^4=\theta(x_3,x_4)x_1=a_{11}=\theta(x_4,x_1)x_3=-a_6\\
& &a_8=\theta(x_2,x_3)x_4=\theta(x_3,x_4)x_2=a_{12}=\theta(x_4,x_2)x_3=-a_{10},
\end{eqnarray*}
so that $\theta=\alpha_1\theta_1+\alpha_2\theta_2+\alpha_3\theta_3+\alpha_4\theta_4$ for $\alpha_1=a_1$, $\alpha_2=a_2$, $\alpha_3=a_4$ and $\alpha_4=a_8$. 

To see that $\theta$ is a 2-cocycle, it suffices to prove that $d\theta (x_1,x_2,x_3)x_4=0$ but, since $x_3,x_4$ are in the center of $\fk$, one has $d\theta(x_1,x_2,x_3)(x_4)=-\theta([x_1,x_2],x_3)x_4=-\theta(x_3,x_3)x_4=0$. 

Finally, to prove that equality (\ref{compJ}) holds, take three different vector $x,y,z$ within the basis $\{x_1,x_2=Jx_1,x_3,x_4=Jx_3\}$. One of the three vectors $x,y,z$ must be the image under $J$ of one or the others. We can consider, without lost of generality, that $y=Jx$ so that we have 
\begin{eqnarray*}\theta(Jx,Jy)z+\theta(Jy,Jz)x+\theta(Jz,Jx)y=\theta(Jx,-x)z+\theta(-x,Jz)x+\theta(Jz,Jx)Jx=\\
\theta(Jx,-x)z+\theta(x,-x)Jz+\theta(Jx,Jx)Jz=\theta(x,Jx)z=
 \theta(x,y)z,\end{eqnarray*}
which finishes the proof.\qed

\medskip

The knowledge of a pseudo-Hermitian quadratic Lie algebra let us also construct many other pseudo-Hermitian quadratic Lie algebras by the following method:

\medskip

\begin{prop}\label{tensor} Let  ${\mathcal A}$ be a real associative commutative algebra  admitting a non-degenerate symmetric bilinear form $B:{\mathcal A}\times{\mathcal A}\to\R$ such that $B(ab,c)=B(b,ac)$ for all $a,b,c\in {\mathcal A}$ and let $(\fg,J,\varphi)$ be a pseudo-Hermitian quadratic Lie algebra. 

The tensor product $\fG=\fg\otimes_\R {\mathcal A}$ with the bracket, complex structure $J'$ and binear form $\phi$ defined by
$$ [x\otimes a,y\otimes b]_\fG=[x,y]\otimes ab,\quad
 J'(x\otimes a)=Jx\otimes a,\quad
\phi(x\otimes a,y\otimes b)=\varphi(x,y)B(a,b),
$$
for every $x,y\in\fg$, $a,b\in{\mathcal A}$, is a pseudo-Hermitian quadratic Lie algebra.
\end{prop}
\proof It was already seen in \cite[Prop. A]{hk} that $(\fG,\phi)$ is a quadratic Lie algebra. Let us prove that $J'$ is a complex structure on $\fG$ and that it is $\phi$-skewsymmetric.

We clearly have $J'(Jx\otimes a)=J^2x\otimes a=-x\otimes a$, for all $x\in\fg,a\in {\mathcal A}$, which shows that $(J')^2=-I$. Further, for every $x,y\in\fg$, $a,b\in{\mathcal A}$ we easily prove that
$
N_{J'}(x\otimes a,y\otimes b)=N_j(x,y)\otimes ab=0$, so that $J'$ is a complex structure and we also have
$$\phi(J'(x\otimes a),y\otimes b)=\varphi(Jx,y)B(a,b)=-\varphi(x,Jy)B(a,b)=-\phi(x\otimes a, J'(y\otimes b)).$$
Therefore, $(\fG,J',\phi)$ is a pseudo-Hermitian quadratic Lie algebra as claimed.\qed

\begin{coro} If $(\fg,J,\varphi)$ is a pseudo-Hermitian quadratic Lie algebra, then the underlying real Lie algebra of its complexification $\fG=(\fg^\C)_\R$ is also a pseudo-Hermitian quadratic Lie algebra with the complex structure $J'$ and quadratic form $\phi$ defined by
$$J'(x+iy)=Jx+iJy,\quad \phi(x+iy,z+it)=\varphi(x,z)-\varphi(y,t),$$
for all $x,y,z,t\in\fg$.
\end{coro}
\proof Notice that we can  naturally identify $\fG=\fg\otimes_\R\C$ by considering $x+iy=x\otimes 1+y\otimes i$ for all $x,y\in\fg$. Put $B(a,b)=\mbox{Re}(ab)$ for $a,b\in\C$. It is easy to see that $B$ is a non-degenerate, symmetric, bilinear form on $\C$ and that $B(ab,c)=B(b,ac)$ holds for all $a,b,c\in\C$. It thus suffices to apply Proposition \ref{tensor} and see that, under the identification $\fG=\fg\otimes_\R\C$, the map $J'$ and the bilinear form $\phi$ are as stated.\qed

\begin{example}{\em When $(\fg,J,\varphi)$ is a pseudo-Hermitian quadratic Lie algebra, one can also use Proposition \ref{tensor} to construct a family of nilpotent pseudo-Hermitian quadratic Lie algebras by appropriate tensorizations of $\fg$. 

Let ${\mathcal A}=\R\mbox{-span}\{a,a^2,\dots,a^k\}$ the nilpotent commutative associative algebra defined by the products $a^ia^j=a^{i+j}$ if $i+j\le k$ and $a^ia^j=0$ whenever $i+j> k$. It suffices to bear in mind that the  form $B:{\mathcal A}\times {\mathcal A}\to\R$ defined by bilinearity from
$B(a^i,a^j)=\delta_{i+j,k+1}$, where $\delta_{r,s}$ stands for the Kronecker symbol,
is a non-degenerate symmetric bilinear form on ${\mathcal A}$ and it verifies
$$B(a^ia^j,a^\ell)=\delta_{i+j+\ell,k+1}=B(a^j,a^ia^\ell)$$
for all $1\le i,j,\ell\le k$. 

According to Proposition \ref{tensor} the (nilpotent) Lie algebra $\fG=\fg\otimes_\R {\mathcal A}$ admits a pseudo-Hermitian quadratic structure.
}
\end{example}

\section{Double extension of pseudo-Hermitian quadratic Lie algebras}

In \cite{MR}, the authors propose a method to construct  $n+2$-dimensional quadratic Lie algebras from a given $n$-dimensional one by performing an appropriate central extension and a semidirect product defined by a certain skewsymmetric derivation. Such a method is known as {\it quadratic double extension by a line}. Explicitly, if $(\fg_0,[\cdot,\cdot]_0,\varphi_0)$ is a quadratic Lie algebra and $D$ is a $\varphi_0$-skewsymmetric derivation of $\fg_0$, then one considers the vector space $\fg=\R z\oplus\fg_0\oplus\R v$ with the bracket $[\cdot,\cdot]$ and the bilinear map $\varphi:\fg\times\fg\to\R$ defined by
\begin{eqnarray*}
& & [z,x]=[z,v]=0,\quad [x,y]=[x,y]_0+\varphi_0(Dx,y)z,\quad [v,x]=Dx\\
& & \varphi(z,v)=1, \quad\varphi(x,y)=\varphi_0(x,y),\quad \varphi(z,\R z\oplus\fg_0)=\varphi(v,\fg_0\oplus\R v)=\{0\}
\end{eqnarray*}
for all $x,y\in\fg_0$. It results that $(\fg,[\cdot,\cdot],\varphi)$ is a quadratic Lie algebra. Further, one has the following result:

\begin{prop}[Medina-Revoy]
Every  non-abelian solvable quadratic Lie algebra $\mathfrak{g}$ with $dim(\mathfrak{g}) = n $ is
the double extension of a solvable quadratic Lie algebra $\mathfrak{g}_0$ of dimension $n - 2$ by a line. As a consequence, every  solvable quadratic  Lie algebra is obtained by a sequence of double
extensions starting from  an abelian Lie algebra.
\end{prop}

In this section we will use quadratic double extensions of pseudo-Hermitian quadratic Lie algebras in order to obtain new ones. Explicitly, we have the following result:

\begin{prop} \label{pR1}
Let $\left( \mathfrak{g}_0 , \left[\cdot, \cdot \right]_0 , J_0 , \varphi_0 \right)$ be a peudo-Hermitian quadratic Lie algebra. Let $D,F$ be two $\varphi_0$-skewsymmetric derivations of $\fg_0$ such that $\left[F+J_0D ,J_0\right]=0 $ and $ \left[  F, D \right]  = \ad_{\mathfrak{g}_0} (s_0)$ for a certain $s_0\in\fg_0$.

Define the vector space $\fg$ as the direct sum of $\fg_0$ and a $4$-dimensional space with generators $\left\lbrace  z, z' , v, v' \right\rbrace  $, namely
$$ \mathfrak{g} = \mathbb{R}z \oplus \mathbb{R}z' \oplus \mathfrak{g}_0 \oplus \mathbb{R}v' \oplus \mathbb{R}v.$$
Let us define on $\fg$ a skewsymmetric product $[\cdot,\cdot]$ and a symmetric bilinear form $\varphi$  by the following nonzero pairings:  
\begin{eqnarray*}
& & [v,v']=s_0,\quad [v,x]=F x - \varphi_0 (s_0 , x)z',\quad [v',x]=Dx + \varphi_0 (s_0 , x)z,\\  & & [x,y]=\left[x, y \right]_0 + \varphi_0 (Dx, y)z' + \varphi_0 (F x, y)z\\
& & \varphi(z, v) = \varphi(z' , v' ) = 1,\quad \varphi(x, y) = \varphi_0 (x, y),
\end{eqnarray*}
for any $x,y\in\fg_0$. Let us also consider $J\in\mathfrak{gl}(\fg)$ defined by
 $$Jz = z',\quad Jz' = -z,\quad  Jv = v',\quad Jv' = -v,\quad Jx = J_0 x,$$ for all $x\in\fg_0$.

The triple $(\fg, [\cdot,\cdot],\varphi)$ is a quadratic Lie algebra, $J$ is a complex structure on $\fg$ and $\varphi$ is actually pseudo-Hermitian.
\end{prop}
\proof Denote by $(\fg_1,\varphi_1)$ the quadratic double extension of $(\fg_0,\varphi_0)$ by means of $D$. One then has $\fg_1=\R z'\oplus\fg_0\oplus\R v'$ with  the nonzero brackets and pairings of $\varphi_1$ given by
$$[x,y]_1=[x,y]_0+\varphi_0(Dx,y)z',\quad [v',x]_1=Dx,\quad \varphi_1(x,y)=\varphi_0(x,y),\quad \varphi_1(z',v')=1,$$
where $x,y\in\fg_0$.
Let $D_1:\fg_1\to\fg_1$ be the linear map defined by 
$$D_1(z')=0,\quad D_1(v')=s_0,,\quad D_1(x)=Fx-\varphi_0(s_0,x)z',$$
for all $x\in\fg_0$. As $F$ is $\varphi_0$-skewsymmetric. for every $\alpha,\alpha',\beta,\beta'\in\R$ and $x,y\in\fg_0$ one has
\begin{eqnarray*}
& & \varphi_1(D_1( \alpha z'+x+\beta v'),\alpha' z'+y+\beta' v')+\varphi_1(\alpha z'+x+\beta v',D_1( \alpha' z'+y+\beta' v'))=\\
& & \varphi_1(Fx-\varphi_0(s_0,x)z'+\beta s_0,y+\beta' v')+\varphi_1(x+\beta v',Fy-\varphi_0(s_0,y)z'+\beta' s_0,)=\\
& & \varphi_0(Fx,y)+\beta\varphi_0(s_0,y)-\beta'\varphi_0(s_0,x)+\varphi_0(x,Fy)+\beta'\varphi_0(x,s_0)-\beta\varphi_0(s_0,y)=0,\end{eqnarray*}
which shows that $D_1$ is skewsymmetric with respect to $\varphi_1$. Let us see that $D_1$ is a derivation of $\fg_1$. One clearly has $D_1[z',X]=[D_1z',X]+[z',D_1X]$ for all $X\in\fg_1$ and for every $x\in\fg_0$ we have that
$$ D_1[v',x]_1-[D_1v',x]_1-[v',D_1x]_1=FDx-\varphi_0(s_0,Dx)z'-[s_0,x]_0-\varphi_0(Ds_0,x)z'-DFx=0$$
because $D$ is $\varphi_0$-skewsymmetric and $[F,D]=\ad_{\fg_0}(s_0)$. Besides, using again this last equality and the facts that $F$ is a skewsymmetric derivation of $(\fg_0,\varphi_0)$ and $\varphi_0$ is quadratic, we  have for all $x,y\in\fg_0$ that
\begin{eqnarray*}
& & D_1[x,y]_1-[D_1x,y]_1-[x,D_1y]_1=\\& & F[x,y]_0-\varphi_0(s_0,[x,y])z'-[Fx,y]_0-\varphi_0(DFx,y)z'-[x,Fy]_0-\varphi_0(Dx,Fy)z'=\\
& & -\varphi_0([s_0,x]y)z'+\varphi_0([F,D]x,y)z'=0.
\end{eqnarray*}
Thus, $D_1$ is a skewsymmetric derivation of $(\fg_1,\varphi_1)$. Now, it is a simple calculation to see that the quadratic double extension of $(\fg_1,\varphi_1)$ by means of $D_1$ is nothing but the pair $(\fg,\varphi)$ of the statement, proving that $(\fg,\varphi)$ is a quadratic Lie algebra. 

Therefore, it only remains to prove that $J$ is a complex structure and that it is $\varphi$-skewsymmetric. It is clear from the definition of $J$ that $J^2=-I$ and we have for $x_k\in\fg_0$ and $\alpha_k,\alpha'_k,\beta_k\beta'_k\in\R$ $(k=1,2)$ that the following holds
\begin{eqnarray*}
& & \varphi(J( \alpha_1z+\alpha_1' z'+x_1+\beta'_1 v'+\beta_1 v),\alpha_2z+\alpha_2' z'+x_2+\beta'_2 v'+\beta_2 v)=\\
& & \varphi( \alpha_1z'-\alpha_1' z+J_0x_1-\beta'_1 v+\beta_1 v',\alpha_2z+\alpha_2' z'+x_2+\beta'_2 v'+\beta_2 v)=\\
& & \alpha_1\beta'_2-\alpha'_1\beta_2-\beta'_1\alpha_2+\beta_1\alpha'_2+\varphi_0(J_0x_1,x_2).
\end{eqnarray*}
Similarly we get 
\begin{eqnarray*}
& &\varphi(\alpha_1z+\alpha_1' z'+x_1+\beta'_1 v'+\beta_1 v,J( \alpha_2z+\alpha_2' z'+x_2+\beta'_2 v'+\beta_2 v))=\\
& & \alpha_2\beta'_1-\alpha'_2\beta_1-\beta'_2\alpha_1+\beta_2\alpha'_1+\varphi_0(J_0x_2,x_1),\end{eqnarray*}
and thus $J$ results to be skewsymmetric with respect to $\varphi$ because $J_0$ is $\varphi_0$-skewsymmetric. In order to see that the $N_J=0$, first notice that one obviously has $N_J(z,\cdot)=N_J(z',\cdot)=0$ because $z,z'=Jz$ are in the center of $\fg$ and that also $N_J(v.v')=0$ because $v'=Jv$. If $x\in\fg_0$ we get
\begin{eqnarray*}
& & N_J(v,x)=\\& & Fx-\varphi_0(s_0,x)z'+J_0Dx+\varphi_0(s_0,x)z'+J_0FJ_0x+\varphi_0(s_0,J_0x)z-DJ_0x-\varphi_0(s_0,J_0x)z=\\
& & (F+J_0D+J_0FJ_0-DJ_0)x=J_0[F+J_0D,J_0]x=0.
\end{eqnarray*}
because $F+J_0D$ commutes with $J_0$ by hypothesis. Notice that, using (\ref{propnijen}), this also implies $N_J(v',x)=N_J(Jv,x)=-JN_J(v,x)=0$.
Further, for $x,y\in\fg_0$ we have
\begin{eqnarray*}
& & [x,y]=[x,y]_0+\varphi_0(Dx,y)z'+\varphi_0(Fx,y)z\\
& & J[Jx,y]=J_0[J_0x,y]_0-\varphi_0(DJ_0x,y)z+\varphi_0(FJ_0x,y)z'\\
& & J[x,Jy]=J_0[x,J_0y]_0+\varphi_0(J_0Dx,y)z-\varphi_0(J_0Fx,y)z'\\
& & [Jx,Jy]=[J_0x,J_0y]_0-\varphi_0(J_0DJ_0x,y)z'-\varphi_0(J_0FJ_0x,y)z,\end{eqnarray*}
from where we immediately get
$$N_J(x,y)=N_{J_0}(x,y)+\varphi_0([F+J_0D,J_0]x,y)z'-\varphi_0(J_0[F+J_0D,J_0]x,y)z=0.$$
Thus, the Nijenhuis tensor vanishes identically and we get that $(\fg,\varphi,J)$ is a pseudo-Hermitian quadratic Lie algebra.\qed

\begin{defi}{\em The pseudo-Hermitian quadratic Lie algebra $(\fg,J,\varphi)$ constructed as in Proposition \ref{pR1} will be called the {\it pseudo-Hermitian 	quadratic double extension of $(\fg_0,J_0,\varphi_0)$ by a plane by means of $(D,F)$ and $s_0$}. For short, we will say that $\fg$ is a {\it pHQ-double extension} of $\fg_0$.}
\end{defi}

\begin{rema}\label{rM1} {\em Let us  remark some interesting facts.
\begin{enumerate}
\item[(1)] Notice that the concept of pHQ-double extension defined above is different from the construction that Andrada, Barberis and Ovando give in \cite[Th. 21]{Andrada}. For a nilpotent Lie algebra $\fg$, their construction is the particular case of ours for which $D=0$. 
\item[(2)] If $\fg$ is a pHQ-double extension of a pseudo-Hermitian quadratic  Lie algebra $\fg_0$ and the signature of the metric in $\fg_0$ is $\mbox{sig}(\varphi_0)=(r,s)$, then for the metric in $\fg$ we have $\mbox{sig}(\varphi)=(r+2,s+2)$.
\item[(3)] If $F,D$ verify $\left[F+J_0D ,J_0\right]=0 $ and $ \left[  F, D \right]  = \ad_{\mathfrak{g}_0} (s_0)$ then for the maps $F_1=D$, $D_1=-F$ one also has $[F_1,D_1]=\ad_{\mathfrak{g}_0} (s_0)$ and
$$[F_1+J_0D_1 ,J_0]=[D-J_0F,J_0]=DJ_0-J_0FJ_0-J_0D-F=J_0[F+J_0D,J_0]=0.$$
If we denote $ \mathfrak{g}_1 = \mathbb{R}z_1 \oplus \mathbb{R}z_1' \oplus \mathfrak{g}_0 \oplus \mathbb{R}v'_1 \oplus \mathbb{R}v_1$ the pHQ-double extension by means of $(D_1,F_1)$ and $\fg$ is as in Proposition \ref{pR1}, one easily sees that the map $\Psi:\fg_1\to\fg$ defined by $\Psi(z_1)=z',\psi(z_1')=-z,\Psi(v_1)=v',\psi(v_1')=-v,\Psi(x)=x$ for all $x\in\fg_0$ is a holomorphic isomorphism that preserves the corresponding metrics.
\item[(4)] Notice that the condition $[ F, D]  = \ad_{\mathfrak{g}_0} (s_0)$ does not completely define $s_0$ since the same identity holds for every $s_0'=s_0+w$ if $w\in\fz(\fg)$. The election of $s_0'$ instead of $s_0$ in the definition of the pHQ-double extension may give rise to nonisomorphic Lie algebras. For instance, when $\fg_0=\R^2$ with a Hermitian metric and we take $D=F=0$, the pHQ-double extension is abelian whenever  $s_0=0$ whereas for $s_0\ne 0$ one gets a non-abelian 3-step nilpotent Lie algebra (see, for instance, Corollary \ref{LH} below).
\end{enumerate}}
\end{rema} 

\begin{prop}\label{pR2} Let $(\fg,J,\varphi)$ be a pseudo-Hermitian 	quadratic Lie algebra of $\dim(\fg)=n>4$ and denote by $\fz(\fg)$ its center. 

If there exists an element $z\in \fz(\fg)\cap [\fg,\fg]$ such that $Jz\in\fz(\fg)$, then there exists an $(n-4)$-dimensional pseudo-Hermitian quadratic Lie algebra $(\fg_0,J_0,\varphi_0)$ such that  $\fg$ is a pHQ-double extension of $\fg$.
\end{prop}
\proof  Since $z$ and $Jz$ are linearly independent, we can find an element $v\in\fg$ such that $\varphi(z,v)=1$ and $\varphi(Jz,v)=0$. Now, as $z\in\fz(\fg)\cap [\fg,\fg]$ implies that $z$ is isotropic,  one can use Medina and Revoy's results in \cite{MR}  to see that there exist a  quadratic Lie algebra $(\fg_1,\varphi_1)$ such that $\fg=\R z\oplus\fg_1\R v$ is the quadratic double extension of $\fg_1$ by means of a certain skewsymmetric derivation $D_1$ of $(\fg_1,\varphi_1)$. Notice that, in that case, $\fg_1$ is the $\varphi$-orthogonal of $\R z\oplus\R v$ and, hence, as $\varphi(Jz,v)=0$ and $\varphi(Jz,z)=0$ because $J$ is skewsymmetric with respect to $\varphi$, we have that $Jz\in\fg_1$. Actually, $Jz$ is an isotropic central element of $\fg_1$ because $\varphi(Jz,Jz)=\varphi(z,z)=0$ and $Jz\in\fz(\fg)$. Since we also have
$$\varphi(Jv,z)=-\varphi(v,Jz)=0,\quad \varphi(Jv,v)=-\varphi(v,Jv)=0,\quad \varphi(Jz,Jv)=\varphi(z,v)=1,$$
we get that $Jv\in\fg_1$ and, using again the results of \cite{MR}, we get that there exists a quadratic Lie algebra $(\fg_0,\varphi_0)$ such that $\fg_1=\R Jz\oplus\fg_0\oplus\R Jv$ is the quadratic double extension of $\fg_0$ by means of a skewsymmetric derivation $D$ of $(\fg_0,\varphi_0)$. If we now denote $z'=Jz,v'=Jv$, we get that $\fg=\mathbb{R}z \oplus \mathbb{R}z' \oplus \mathfrak{g}_0 \oplus \mathbb{R}v' \oplus \mathbb{R}v$, where $\fg_0$ is the $\varphi$-orthogonal to $\R\mbox{-span}\{z,z',v,v'\}$. Besides, as $\R\mbox{-span}\{z,z',v,v'\}$ is clearly $J$-invariant.

From the construction of quadratic double extensions we have that the only {\it a priori} nonzero pairings for $\varphi$ are
$$\varphi(z,v)=1,\quad \varphi(z',v')=\varphi_1(z',v')=1,\varphi(x,y)=\varphi_1(x,y)=\varphi_0(x,y), \quad x,y\in\fg_0,$$ and the bracket in $\fg$ is defined by $z,z'\in\fz(\fg)$ and
\begin{eqnarray*}
& & [x,y]_1=[x,y]_0+\varphi_0(Dx,y)z',\quad [v',x]_1=Dx\\
& & [x,y]=[x,y]_1+\varphi_1(D_1x,y)z,\quad [v',x]=[v',x]_1+\varphi_1(D_1v',x)z,\\ & & [v,x]=D_1x,\quad [v,v']=D_1v',
\end{eqnarray*}
for all $x,y\in\fg_0$. Notice that $D_1z'=0$ because, being $D_1$ skewsymmetric with respect to $\varphi_1$, we have $\varphi_1(D_1z',z')=0$ and for $x\in\fg_0$ and $\lambda\in\R$ we get
$$\varphi_1(x+\lambda v',D_1z')=-\varphi_1(D_1x+\lambda D_1v',z')=-\varphi([v,x+\lambda v'],z')=\varphi(x+\lambda v',[v,z'])=0.$$
This implies that $\varphi_1(D_1\fg_1,z')=\varphi_1(\fg_1,D_1z')=\{0\}$, showing that $D_1\fg_1\subset \R z'\oplus\fg_0$. But $\varphi_1(D_1v',v')=0$ because $D_1$ is skewsymmetric and therefore we have:
$$D_1v'\in\fg_0,\quad D_1x=Fx+\alpha(x)z',\quad x\in\fg_0$$
for a certain $\alpha\in\fg_0^*$ and a a certain linear endomorphism $F$ of $\fg_0$. But  then the skewsymmetry of $D_1$ gives that $F$ is $\varphi_0$-skewsymmetric and that
$$\alpha(x) =\varphi_1(D_1x,v')=-\varphi_1(x,D_1v').$$ 
The fact that $D_1$ is a derivation of $\fg_1$ gives
$$F[x,y]_0-\varphi_0([x,y]_0,s_0)z'=[Fx,y]_0+\varphi_0(DFx,y)z'+[x,Fy]_0+\varphi_0(Dx,Fy)z',$$
from where we derive that $F$ is a derivation of $\fg_0$ and that $[F,D]=\mbox{ad}_{\fg_0}(s_0)$ because
$$\varphi_0([s_0,x],y)=\varphi(s_0,[x,y])=-\varphi_0(DFx,y)-\varphi_0(Dx,Fy)=\varphi([F,D]x,y),$$
for all $x,y\in\fg_0$.
Summing up, if we denote $s_0=D_1v'$ we get that the bracket in $\fg$ is given by 
\begin{eqnarray*}
& &[x,y]=[x,y]_0+\varphi_0(Dx,y)z'+\varphi_0(Fx,y)z,\quad [v',x]=Dx+\varphi_1(s_0,x)z,\\ & & [v,x]=Fx-\varphi_1(s_0,x)z',\quad [v,v']=s_0.\end{eqnarray*}

Thus, it only remains to prove that $J_0$ is a $\varphi_0$-skewsymmetric complex structure on $\fg_0$ and that $[F+J_0D,J_0]=0$. Uding the brackets above and the fact that $F$ is skewsymmetric for $ \varphi_0$, we have for $x,y\in\fg_0$ that
\begin{eqnarray*}
& & [x,y]=[x,y]_0+\varphi_0(Dx,y)z'+\varphi_0(Fx,y)z\\
& & J[Jx,y]=J_0[J_0x,y]_0-\varphi_0(DJ_0x,y)z+\varphi_0(FJ_0x,y)z'\\
& & J[x,Jy]=J_0[x,J_0y]_0+\varphi_0(J_0Dx,y)z-\varphi_0(J_0Fx,y)z'\\
& & [Jx,Jy]=[J_0x,J_0y]_0-\varphi_0(J_0DJ_0x,y)z'-\varphi_0(J_0FJ_0x,y)z\end{eqnarray*}
So, $N_J(x,y)=0$ is equivalent to $N_{J_0}(x,y)=0$ and
$$0=D+FJ_0-J_0F+J_0DJ_0=[F+J_0D,J_0],\quad 0=F-DJ_0+J_0D+J_0FJ_0=J_0[F+J_0D,J_0].$$
Since one clearly has $J_0^2x=J^2x=-x$ and 
$\varphi_0(J_0x,y)=\varphi(Jx,y)=-\varphi(x,Jy)=-\varphi_0(x,J_0y),$
for every elements $x,y\in\fg_0$, we can conclude that $(\fg,J,\varphi)$ is the pHQ-double extension of $(\fg_0,J_0,\varphi_0)$ by means of $(D,F)$,\qed

\bigskip

We can now state our main result. Its proof is heavily based on the result of Salamon \cite[Corol. 1.4]{salamon} (see also \cite[Prop. 3.2]{Million}) assuring that for a nilpotent Lie algebra with a complex structure $(\fg, J)$ one always has $[\fg,\fg]+J[\fg,\fg]\ne\fg.$

\begin{theor}\label{mainth} Let $(\fg,J,\varphi)$ be a nilpotent pseudo-Hermitian quadratic Lie algebra ans suppose that $\dim(\fg)> 4$. The following hold: 

\begin{enumerate}
\item[{\rm (1)}] Either $\fg$ is a pHQ-double extension of an $(n-4)$-dimensional pseudo-Hermitian quadratic Lie algebra or it is the orthogonal direct sum of a $(n-2)$-dimensional nilpotent pseudo-Hermitian quadratic Lie algebra and $\R^2$ endowed with a definite pseudo-Hermitian metric.
\item[{\rm (2)}] If $\fg$ is not abelian, then it is a pHQ-double extension of an $(n-4)$-dimensional pseudo-Hermitian quadratic Lie algebra.

As a consequence, every non-abelian pseudo-Hermitian quadratic Lie algebra is obtained from a series of pHQ-double extensions starting from an abelian Lie algebra.
\end{enumerate}
\end{theor}
\proof As proved by Salamon in \cite{salamon}, we have that $[\fg,\fg]+J[\fg,\fg]\ne\fg$ and therefore its orthogonal spaces with respect of $\varphi$ verify $[\fg,\fg]^\bot\cap (J[\fg,\fg])^\bot\ne\{0\}$. Since $[\fg,\fg]^\bot=\fz(\fg)$ and $(J[\fg,\fg])^\bot=J([\fg,\fg]^\bot)=J\fz(\fg)$, we get that $\fz(\fg)\cap J\fz(\fg)\ne\{0\}.$ Let us take a nonzero $z\in \fz(\fg)\cap J\fz(\fg)$. Obviously, we have $Jz\in\fz(\fg)$. If $z$ is isotropic, this means $z\in \fz(\fg)\cap [\fg,\fg]$ and then Proposition \ref{pR2} guarantees that $\fg$ is a pHQ-double extension. Otherwise, $\varphi(z,z)=\varphi(Jz,Jz)\ne 0$ and $\varphi(z,Jz)=0$ so that $\fI=\R\mbox{-span}\{z,Jz\}$ is a non-degenerate  $J$-invariant ideal of $\fg$ and the restriction of $\varphi$ to $\fI\times \fI$ is definite (positive or negative). Further, $\fI\subset \fz(\fg)$ clearly implies $[\fg,\fg]=\fz(\fg)^\bot\subset \fI^\bot$, which proves that $\fI^\bot$ is also an ideal of $\fg$ and, since $J$ is skewsymmetric with respect to $\varphi$, it is $J$-invariant. Thus $\fg=\fI^\bot\oplus\fI=\fI^\bot\oplus\R^2$, which proves (1).

In order to prove (2), notice that if $\fg$ is non-abelian, we can apply recursively (1) to obtain an orthogonal decomposition of Lie algebras $\fg=\fh\oplus (\R^2)^k$ for certain $k\ge 0$ where $\fh$ is the pHQ-double extension of a triple $(\fh_0,J_0,\varphi_0)$ by means of certain $(F_\fh,D_\fh,s_\fh)$. A simple calculation shows that then $\fg$ is the pHQ-double extension of $\fg_0=\fh_0\oplus (\R^2)^k$ by means of the triple $(F,D,s_0)$ defined by $F(h,t)=F_\fh(h),D(h,t)=D_\fh(h)$ for all $h\in\fh_0,t\in (\R^2)^k$ and $s_0=(s_\fh,0)$.

The second part of (2) follows by repeating the same argument to $\fg_0$ and so on. Obviously, the process stops when one reaches an abelian reduced Lie algebra. \qed

\begin{rema}{\em Notice that, if $\fg$ is a pHQ-double extension of another pseudo-Hermitian quadratic Lie algebra $\fg_0$ and $\fg$ is nilpotent then, as stated in the theorem, $\fg_0$ must be nilpotent but also the derivations $D,F$ defining the double extension should be nilpotent. Actually, if $\fg=\mathbb{R}z \oplus \mathbb{R}z' \oplus \mathfrak{g}_0 \oplus \mathbb{R}v' \oplus \mathbb{R}v$ is $k$-step nilpotent then $\ad (v)^k=\ad (v')^k=0$ and $\ad(x)^k=0$ for all $x\in\fg$ but then, for each $y\in\fg_0$, we have
\begin{eqnarray*}
& & 0=\ad(v)^k(y)=F^ky-\varphi_0(s_0,F^{k-1}y)z'\\
& & 0=\ad(v')^k(y)=D^ky+\varphi_0(s_0,D^{k-1}y)z'\\
& & 0=\ad(x)^k(y)=\ad_{\fg_0}(x)^ky+\varphi_0(Dx,\ad_{\fg_0}(x)^{k-1}y)z'++\varphi_0(Fx,\ad_{\fg_0}(x)^{k-1}y)z.\end{eqnarray*}
The two first equalities show that $F$ and $D$ are nilpotent and the third one that $\fg_0$ is a nilpotent Lie algebra.
}\end{rema}

Using the theorem and the remark above we  can give a complete description of Lorentz-Hermitian quadratic Lie algebras.

\begin{coro}\label{LH} A non-abelian nilpotent pseudo-Hermitian quadratic Lie algebra $(\fg,J,\varphi)$ is Lorentz-Hermitian if and only if it is the orthogonal sum of an abelian Hermitian Lie algebra and the $6$-dimensional Lorentz-Hermtian Lie algebra ${\mathcal L}=\R\mbox{\rm{-span}}\{x_1,Jx_1,x_2,Jx_2,x_3,Jx_3\}$ with nonzero brackets and $\varphi$-pairings given by
\begin{eqnarray*}
& & [x_1,Jx_1]=x_2, \quad [x_1,x_2]=- Jx_3,\quad [Jx_1,x_2]=x_3\\
& & \varphi(x_1,x_3)=\varphi(Jx_1,Jx_3)=1,\quad \varphi(x_2,x_2)=\varphi(Jx_2,Jx_2)=1.
\end{eqnarray*}
\end{coro}
\proof Since up to dimension $4$ every nilpotent quadratic Lie algebra is abelian \cite{FS,IK}, we can consider that $\dim(\fg)=n\ge 6$. According to Theorem \ref{mainth}, $\fg$ is  a pHQ-double extension of a $(n-4)$-dimensional pseudo-Hermitian quadratic Lie algebra. 

Recalling Remark \ref{rM1} (2), the corresponding Lie algebra $\fg_0$ must have signature $(0,n-4)$. But the only nilpotent quadratic Lie algebras with definite positive metric are the abelian ones, so that $\fg_0=\R^{n-4}$ with an arbitrary complex structure and a Hermitian metric. It is well known that a nilpotent linear map which is skew-symmetric with respect to a definite scalar product must be identically zero. So, the derivations $D$ and $F$ defining the pHQ-double extension are both null. This means that the only nonzero brackets in $\fg=\mathbb{R}z \oplus \mathbb{R}z' \oplus \mathfrak{g}_0 \oplus \mathbb{R}v' \oplus \mathbb{R}v$ are
$$ [v,v']=s_0,\quad [v,x]=-\varphi(s_0,x)z',\quad [v',x]=\varphi(s_0,x)z.$$
Put $n=2k$, $\alpha=\varphi(s_0,s_0)$ and choose a basis $\{Js_0,s_1,Js_1,\dots,s_{k-3},Js_{k-3}\}$ of $s_0^\bot$ in $\fg_0$. Then, the nonzero brackets in $\fg$ are
$$[v,v']=s_0,\quad [v,s_0]=-\alpha z',\quad [v',s_0]=\alpha z.$$

Notice that $\alpha>0$ because $\varphi_0$ is positive definite. Put  $x_1=\alpha^{-1/4}v$, $x_2=\alpha^{-1/2}s_0$ and $x_3=\alpha^{1/4}z$. Since $v'=Jv$ and $z'=Jz$, we get that the non-null brackets are
$$[x_1,Jx_1]= \alpha^{-1/2}s_0=x_2,\,\, [x_1,x_2]=-\alpha^{-1/4}\alpha^{-1/2}\alpha Jz=-Jx_3,\,\, 
[Jx_1,x_2]=\alpha^{-1/4}\alpha^{-1/2}\alpha z=x_3.$$
Further, in ${\mathcal L}=\R\mbox{\rm{-span}}\{x_1,Jx_1,x_2,Jx_2,x_3,Jx_3\}$ the metric is given by
$$\varphi(Jx_1,Jx_3)=\varphi(x_1,x_3)=\alpha^{-1/4}\alpha^{1/4}\varphi(v,z)=1,\,\, \varphi(Jx_2,Jx_2)=\varphi(x_2,x_2)=\alpha^{-1}\varphi(s_0,s_0)=1.$$
Now, putting  ${\mathcal A}= \R\mbox{-span}\{s_1,Js_1,\dots,s_{k-3},Js_{k-3}\}$, one immediately gets that $\fg$ is the orthogonal sum $\fg={\mathcal L}\oplus {\mathcal A}$ where ${\mathcal A}$ is an abelian Lie algebra with Hermitian metric.
\qed

\begin{rema}{\em The Lie algebra ${\mathcal L}$ is isomorphic to ${\mathcal W}_3\oplus \R$ with the notation of  \cite{FS} and it is the Lie algebra denoted $(0,0,0,12,14,24)$ in the list of $6$-dimensional nilpotent Lie algebras admitting complex structures given by Salamon in \cite{salamon}. In \cite{Andrada}, it is denoted by $L_3(1,2)\times\R.$}
\end{rema}
 
\section{Classification of pseudo-Hermitian quadratic nilpotent Lie algebras up to dimension $8$}\label{classif}

We will now give a complete classification of the nilpotent Lie algebras which admit a pseudo-Hermitian quadratic metric. Notice that our classification actually done up to  holomorphic isometric isomorphism but in the cases of abelian Lie algebras we will not list all the obvious distinct signatures. 

The classification up to dimension $6$ is now almost obvious. With the same notation as in Corollary \ref{LH}, one has

\begin{prop} A nilpotent Lie algebra $\fg$ with $\dim(\fg)\le 6$ admits a pseudo-Hermitian quadratic structure if and only if it is one of the abelian Lie algebras $\R^2$, $\R^4$, $\R^6$ with a metric of signature $(2r,2p)$ or ${\mathcal L}$ with either the metric of Corollary \ref{LH} or its opposite.
\end{prop}
\proof Recall that up to dimension $4$ all the nilpotent quadratic Lie algebras are abelian \cite{FS} and a pseudo-Hermitian quadratic Lie algebra has even dimension. So up to dimension $4$ we only have the Lie algebras $\R^2$, $\R^4$. In the $6$-dimensional case, by Theorem \ref{mainth} we have that a pseudo-Hermitian quadratic Lie algebra $\fg$ is either the sum of $\R^2$ and a $4$-dimensional pseudo-Hermitian quadratic Lie algebra or a pHQ-double extension of $\R^2$ with a definite pseudo-Hermitian metric.  In the first case we get $\fg=\R^6$ and in the second (changing the sign to the metric on $\fg$ if necessary) the Lie algebra ${\mathcal L}$ as in the proof of Corollary \ref{LH} if one chooses $s_0\ne 0$ and $\R^6$ when $s_0=0$.\qed

\bigskip

Let us then study the case $\dim(\fg)=8$. Using Theorem \ref{mainth} again and the classification in dimension $6$ we get that if  $\fg$ is not a pHQ-double extension of $\R^4$ then it must be isomorphic to $\R^8$ or to ${\mathcal L}\oplus\R^2$. So, we can consider that $\fg$ is a pHQ-double extension of $\R^4$. Notice that for a definite metric on $\R^4$ we also get either $\R^8$ or the orthogonal sum ${\mathcal L}\oplus\R^2$ (with the metric given in Corollary \ref{LH} or its opposite) so that we can consider that $\fg$ is a pHQ-double extension of $\R^4$ with a neutral pseudo-Hermitian metric (denoted, as usual, $\R^{2,2}$).

\medskip

Let us start with the case in which the pHQ-doble extension is done with two null derivations.

\begin{prop} Let $(\fg,J,\varphi)$ be the  pHQ-doble extension of $\R^{2,2}$ by means of $D=F=0$ and $s_0\in\R^4$. Then it is equivalent to either $\R^{4,4}$ or $T^*_0\fk$ or else $T^*_{\pm\theta_1}\fk$, where $\fk$ and $\theta_1$ are as in Proposition \ref{kodaira} .
\end{prop}
\proof From the definition of the pHQ-double extension we have that the non-zero brackets are
$$[v,Jv]=s_0,\quad [v,x]=-\varphi(s_0,x)Jz,\quad [Jv,x]=\varphi_0(s_0,x)z.$$
If $s_0=0$ we obviously get the abelian Lie algebra $\R^8$ with a neutral metric so that we can consider $s_0\ne 0$.

Suppose that $\varphi_0(s_0,s_0)=0$. We may then find an isotropic vector $t_0\in\R^4$ such that $\varphi_0(s_0,t_0)=1,\varphi_0(Js_0,t_0)=0$. Our non-trivial brackets in the basis $ \{v,Jv,z,Jz,s_0,Js_0,t_0,Jt_0\}$ are then
$$[v,Jv]=s_0,\quad [v,t_0]=-Jz,\quad [Jv,t_0]=z.$$
If we define
$$x_1=v,\quad x_2=Jv,\quad x_3=s_0,\quad x_4=Js_0,\quad x_1^*=z,\quad x_2^*=Jz,\quad x_3^*=t_0,\quad x_4^*=Jt_0,$$
we immediately get that $\fg$ is the Lie algebra $T^*_0\fk$ with the pseudo-Hermitian structure of Proposition \ref{kodaira}.

On the other hand, if $\varphi(s_0,s_0)=\alpha\ne 0$, since the metric in $\R^4$ is neutral we can choose $t_0\in (\R s_0\oplus\R Js_0)^\bot$ such that $\varphi_0(t_0,t_0)=-\varphi_0(s_0,s_0)$ so that in the corresponding basis we get the non-zero brackets
\begin{equation}\label{eqL} [v,Jv]=s_0,\quad [v,s_0]=-\alpha Jz,\quad [Jv,s_0]=\alpha z.\end{equation}
Choose $\lambda\in\R$ such that $\alpha=\epsilon 2\lambda^4$ where $\epsilon=\pm1$ is the sign of $\alpha$ and define
\begin{eqnarray*}
& & x_1=\lambda^{-1}v,\quad  x_2=\lambda^{-1}Jv,\quad x_3=(2\lambda^2)^{-1}(s_0+t_0),\quad x_4=(2\lambda^2)^{-1}(Js_0+Jt_0),\\
& &  x_1^*=\lambda z,\quad x_2^*=\lambda Jz,\quad x_3^*=\epsilon(2\lambda^2)^{-1}(s_0-t_0),\quad x_4^*=\epsilon(2\lambda^2)^{-1}(Js_0-Jt_0).\end{eqnarray*}
We then have that $x_1^*,x_2^*,x_4,x_4^*\in\fz(\fg)$ and a simple calculation yields
$$[x_1,x_2]=x_3+\epsilon x_3^*,\quad [x_2,x_3]=\epsilon x_1^*,\quad [x_3,x_1]=\epsilon x_2^*,\quad [x_3^*,x_1]= x_2^*,\quad [x_3^*,x_2]= -x_1^*,$$
which shows that the bracket of $\fg$ is that of $T^*_{\epsilon \theta_1}\fk$. Since the complex structure is obviously the same and we also have that the non-zero products for $\varphi$ are $\varphi(x_i,x_i^*)=1$, we are done.\qed

\begin{rema}{\em Notice that from the identity (\ref{eqL}) one easily sees that $T^*_{\pm\theta_1}\fk$ is holomorphically isomorphic to ${\mathcal L}\oplus\R^2$. Further, if one chooses on ${\mathcal L}$ the Lorentz-Hermitian metric of Corollary \ref{LH} one has that $T^*_{\theta_1}\fk$ is equivalent to ${\mathcal L}\oplus\R^{2,0}$ and $T^*_{-\theta_1}\fk$ is the same Lie algebra with the opposite metric.  Besides, $T^*_{\theta_1}\fk$ and $T^*_{-\theta_1}\fk$ cannot be equivalent since the restriction of the metric to the derived ideal has signature $(0,1)$ for $T^*_{\theta_1}\fk$ and signature $(1,0)$ for $T^*_{-\theta_1}\fk$.}\end{rema}

In order to classify the cases with $(D,F)\ne(0,0)$, we will always consider that $F\ne 0$. Notice that, according to Remark \ref{rM1} (3), there is no lost of generality with that assumption. We first need the following technical lemma:
  
\begin{lemma}\label{algeb} In $\R^4$ with an arbitrary complex structure $J$ and a  neutral pseudo-Hermitian  metric $\varphi$, let us consider  $F,D\in {\mathfrak{gl}}(\R^4)$ two nilpotent $\varphi$-skewsymmetric maps such that $F\ne 0$  and $[F,D]=[F+JD,J]=0.$ 

If we denote by $\mbox{\rm Ker}(F), \mbox{\rm Im}(F)$ respectively the kernel and the image of $F$ then we have
$$\mbox{\rm Ker}(F)=J(\mbox{\rm Ker}(F))=\mbox{\rm Im}(F).$$
\end{lemma}
\proof Notice that $F$ is skewsymmetric with respect to $\varphi$ and hence must have even rank. Since $F\ne 0$ and $F$ is singular (because it is nilpotent), we must have $\mbox{rank}(F)=2$. This obviously implies that $\mbox{Ker}(F)$ is  $2$-dimensional. Also recall that, if $F^*$ denotes the $\varphi$-adjoint for $F$, since $F$ is $\varphi$-skewsymmetric, we have $\mbox{Im}(F)=\mbox{Ker}(F^*)^\bot=\mbox{Ker}(F)^\bot$.

Let us first suppose that $\mbox{Ker}(F)\cap J(\mbox{Ker}(F))=\{0\}$. Then $\R^4=\mbox{Ker}(F)\oplus J(\mbox{Ker}(F))$. Since $F$ is nilpotent, we may find a non-zero element $t_1\in\mbox{Ker}(F)\cap\mbox{Im}(F)$ and consider another $t_2\in \mbox{Ker}(F)$ linearly independent with $t_1$. We then have that $\{t_1,t_2,Jt_1,Jt_2\}$ is a basis of $\R^4$. Let us denote by ${\mathcal F},{\mathcal D}, {\mathcal J}$ and ${\mathcal G}$ respectively the corresponding matrices of $F,D,J$ and $\varphi$ with respect to such a basis. As $t_1\in\mbox{Im}(F)$ it must be orthogonal to $\mbox{Ker}(F)$ and, therefore, we have
$$\varphi(t_1,t_1)=\varphi(t_1,t_2)=\varphi(Jt_1,Jt_1)=\varphi(Jt_1,Jt_2)=\varphi(t_1,Jt_1)=0.$$
Besides, there exist $\alpha,\beta\in\R$ with $\beta\ne 0$ such that 
$$\varphi(t_2,t_2)=\varphi(Jt_2,Jt_2)=\alpha,\quad \varphi(t_1,Jt_2)=-\varphi(Jt_1,t_2)=\beta,$$
so that the matrices ${\mathcal G}$ and ${\mathcal J}$ are given by
$${\mathcal G}=\left(\begin{array}{rrrr}0 &0 &0&\beta\\0&\alpha&-\beta&0\\0&-\beta&0&0\\\beta&0&0&\alpha\end{array}\right),\quad 
{\mathcal J}=\left(\begin{array}{rrrr}0 &0 &-1&0\\0&0&0&-1\\1&0&0&0\\0&1&0&0\end{array}\right).$$
On the other hand, there exists two $2\times 2$ matrices $S_1$, $S_2$ such that ${\mathcal F}$ has the following block decomposition
$${\mathcal F}=\left(\begin{array}{c|c} 0&S_1\\ \hline 0&S_2\end{array}\right),$$
but one easily sees that the condition $F+F^*=0$ or, equivalently, ${\mathcal G}{\mathcal F}+{\mathcal F}^t{\mathcal G}=0$, implies that there exists $a\in\R$ such that 
$$S_1=aI,\quad S_2=\left(\begin{array}{cc}0&a\alpha/\beta\\0&0\end{array}\right).$$
Now, since $[F,D]=0$ one has that $D(\mbox{Ker}(F))\subset \mbox{Ker}(F)$. So, the matrix of $D$ with respect to the given basis must be of the form
$${\mathcal D}=\left(\begin{array}{c|c} A_1&A_2\\ \hline 0&A_3\end{array}\right)$$ for some $2\times 2$ matrices $A_1,A_2,A_3.$
The condition $[F+JD,J]=0$ leads to
$A_1=A_3-aI$ and $A_2=-S_2.$ But $F$ is nilpotent and this implies that both $A_3$ and $A_1=A_3-aI$ must be nilpotent, hence traceless, and this would imply $a=0$, a contradiction with the assumption $F\ne 0$.

As a consequence, we finally have $\mbox{Ker}(F)\cap J(\mbox{Ker}(F))\ne\{0\}$ so that we can find a vector $u\in \mbox{Ker}(F)$ with $Ju\in \mbox{Ker}(F)$ and this means $\mbox{Ker}(F)=\R\mbox{-span}\{u,Ju\}=J(\mbox{Ker}(F))$. If $\varphi(u,u)=\alpha\ne 0$, we would have that $\varphi$ is definite on $\mbox{Ker}(F)$ but this is impossible  since $\mbox{Ker}(F)\cap\mbox{Im}(F)\ne\{0\}$ by the nilpotency of $F$ and, therefore, there exist an isotropic vector in $\mbox{Ker}(F)$. Thus, $\varphi(u,u)=0$ and $\mbox{Ker}(F)=\mbox{Ker}(F)^\bot=\mbox{Im}(F)$.\qed

\begin{rema}{\em Recall that on a pseudo-Hermitian Lie algebra $(\fg, J,\varphi)$ if $u_1\in\fg$ is isotropic then the subspace $\R u_1\oplus\R Ju_1$ is totally isotropic. Since $u_1^\bot\ne (Ju_1)^\bot$, we may find $v\in (Ju_1)^\bot$ such that $\varphi(u_1,v)=\lambda\ne 0$ so that $u_2=-\frac{\varphi(v,v)}{2\lambda^2}u_1+\frac{1}{\lambda}v$ verifies $\varphi(u_2,u_2)=0$,
$\varphi(u_1,u_2)=1$ and $\varphi(Ju_1,u_2)=0$.
}\end{rema}

\begin{lemma}\label{algeb2} Let us consider $\R^4$ with an arbitrary complex structure $J$ and a  neutral pseudo-Hermitian  metric $\varphi$. Let $F,D$ be two nilpotent $\varphi$-skewsymmetric maps such that $F\ne 0$  and $[F,D]=[F+JD,J]=0$ and consider a non-zero vector $u_1\in\mbox{\rm Ker}(F)$ and a non-zero isotropic $u_2\in\R^4$ such that $\varphi(u_1,u_2)=1$, $\varphi(Ju_1,u_2)=0$. 

The set $\{u_1,Ju_1,u_2,Ju_2\}$ is a basis of $\R^4$, the metric $\varphi$ is given by the only non-zero pairings
$\varphi(u_1,u_2)=\varphi(Ju_1,Ju_2)=1,$ and
there exist $a,b\in\R$ with $a\ne 0$ such that
$F(u_2)=aJu_1$, $F(Ju_2)=-au_1$, $D(u_2)=bJu_1$, $D(Ju_2)=-bu_1$ and, besides, $D(u_1)=D(Ju_1)=0$. 
\end{lemma}
\proof  The facts that $\{u_1,Ju_1,u_2,Ju_2\}$ is a basis of $\R^4$ and that $\varphi(u_1,u_2)=\varphi(Ju_1,Ju_2)=1$ are quite obvious so that we will only see the last part of the statement.

Notice that, according to the lemma above, $\mbox{Ker}(F)=\R\mbox{-span}\{u_1,Ju_1\}=\mbox{Im}(F)$ so that $F(u_2)=\alpha_1u_1+\alpha_2Ju_1$, $F(Ju_2)\beta_1u_1+\beta_2Ju_1$ but since $F$ is skewsymmetric with respect to $\varphi$ we have
$\alpha_1=\varphi(F(u_2),u_2)=0$, $\beta_2=\varphi(F(Ju_2),Ju_2)=0$ and$$\alpha_2=\varphi(F(u_2), Ju_2)=-\varphi(u_2, F(Ju_2))=-\beta_1,$$
so that putting $a=\alpha_2$ we get $F(u_2)=aJu_1$, $F(Ju_2)=-au_1$. Remark that we then have $[F,J]=0$.

The condition $[F,D]=0$ shows that $D(\mbox{Ker}(F))\subset \mbox{Ker}(F)$ so that the matrix of $F$ with respect to the basis $\{u_1,Ju_1,u_2,Ju_2\}$ is of the form 
$${\mathcal D}=\left(\begin{array}{c|c} A_1&A_2\\ \hline 0&A_3\end{array}\right)$$ for some $2\times 2$ matrices $A_i$. As we had obtained $[F,J]=0$,  the condition $[F+JD,J]=0$ implies $[D,J]=0$ and a simple computation shows that each of the matrices $A_i$ must be of the form
$$A_i= \left(\begin{array}{rr}a_i&-b_i\\b_i&a_1\end{array}\right),$$
for certain $a_i,b_i\in\R$. The nilpotency of $D$ now implies that $a_1=b_1=a_3=b_3=0$, showing that $\mbox{Ker}(F)\subset \mbox{Ker}(D)$ and since $D$ is $\varphi$-skewsymmetric we have
$a_2=\varphi(F(u_2),u_2)=0$, so that the election $b=b_2$ gives the desired result.\qed 

\begin{prop} Let $(\fg,J,\varphi)$ by a pHQ-double extension of $\R^{2,2}$ by means of a pair $(F,D)\ne (0,0)$ and a certain $s_0\in\R^4$ and let $\fk$  and $\theta_i$  $(1\le i\le 4)$ be as in Proposition \ref{kodaira}. Then $(\fg,J,\varphi)$ is equivalent to   either $T^*_0\fk$ or $T^*_{\theta_3}\fk$.
\end{prop}
\proof 
As seen in Remark \ref{rM1} (3), we can consider that $F\ne 0$. Lemma \ref{algeb} clearly implies that $\mbox{Ker}(F)$ is 2-dimensional and totally isotropic. Let us consider a 2-dimensional totally isotropic complement $W$ of $\mbox{Ker}(F)$ in $\R^4$ and put $s_0=t_0+w_0$ with $t_0\in\mbox{Ker}(F),w_0\in W$. We will now distiguish  several cases  depending on the vanishing or not of $s_0$ and $w_0$.

If $s_0=0$, let us take an arbitary $u_1\in\mbox{Ker}(F)$ and choose  $u_2\in W$ such that 
$\varphi_0(u_1,u_2)=1$ and $\varphi_0(Ju_1,u_2)=0$. Using Lemma \ref{algeb2} and the bracket given in Proposition \ref{pR2} we get that there exists $a,b\in\R$, $a\ne 0$ such that the brackets in $\fg$ are given by
\begin{eqnarray*}
& & [v,Jv]=0,\quad [v,u_2]=aJu_1,\quad [v,Ju_2]=-au_1\\
& & [Jv,u_2]=bJu_1,\quad [Jv,Ju_2]=-bu_1,\quad [u_2,Ju_2]=az+bJz.
\end{eqnarray*}
Recall that the quadratic form on $\fg$ is given by the non-zero products
$$\varphi(v,z)=\varphi(Jv,Jz)=\varphi(u_1,u_2)=\varphi(Ju_1,Ju_2)=1.$$
Let us now define the vectors
\begin{eqnarray*}
& &x_1=u_2,\quad x_3=az+bJz,\quad x_1^*=u_1,\quad x_3^*=(a^2+b^2)^{-1}(av+bJv),\\ & &  x_2=Jx_1,\quad x_4=Jx_3,\quad x_2^*=Jx_1^*,\quad x_4^*=Jx_3^*.\end{eqnarray*}
It is obvious that $x_1^*,x_2^*,x_3,x_4$ are in the center of $\fg$ and an easy calculation shows that also $x_4^*$ is central and then the only non-trivial brackets are
$$[x_1,x_2]=x_3,\quad [x_3^*,x_1]=x_2^*,\quad[x_3^*,x_2]=-x_1^*.$$
Since one easily verifies that $\varphi(x_i,x_j^*)=\delta_{ij}$ and $\varphi(x_i,x_j)=\varphi(x_i^*,x_j^*)=0$ for all $1\le i,j\le 4$, we see that $\fg=T^*_0\fk$ with the bracket and complex structure defined in Proposition \ref{kodaira}.

If $s_0\ne 0$ but $w_0=0$ we can choose $u_1=t_0\in\mbox{Ker}(F)$ and find a $u_2\in W$ verifying 
$\varphi_0(u_1,u_2)=1$ and $\varphi_0(Ju_1,u_2)=0$.  Combining Lemma \ref{algeb2} and  Proposition \ref{pR2} one obtains  $a,b\in\R$, $a\ne 0$ such that the brackets  are
\begin{eqnarray*}
& & [v,Jv]=u_1,\quad [v,u_2]=aJu_1-Jz,\quad [v,Ju_2]=-au_1\\
& & [Jv,u_2]=bJu_1+z,\quad [Jv,Ju_2]=-bu_1,\quad [u_2,Ju_2]=az+bJz.
\end{eqnarray*}
Consider the following vectors in $\R^4$:
\begin{eqnarray*}
& & x_1=u_2-av-bJv,\quad x_3=2(a^2+b^2)u_1+2az+2bJz,\\
& & x_1^*=\frac{1}{2a^2+2b^2}((a^2+b^2)u_1-az-bJz)\,\quad x_3^*=\frac{1}{4a^2+4b^2}(u_2+av+bJv)\\
& &x_2=Jx_1,\quad x_4=Jx_3,\quad x_2^*=Jx_1^*,\quad x_4^*=Jx_3^*.
\end{eqnarray*}
As in the previous case, it is clear that $x_1^*,x_2^*,x_3,x_4\in\fz(\fg)$. A straightforward computation gives that $x_4^*\in\fz(\fg)$ and also
\begin{eqnarray*}
& & [x_1,x_2]=x_3,\quad [x_3^*,x_1]=x_2^*,\quad[x_3^*,x_2]=-x_1^*\\
& & \varphi(x_i,x_j^*)=\delta_{ij},\quad \varphi(x_i,x_j)=\varphi(x_i^*,x_j^*)=0,\quad 1\le i,j\le 4,
\end{eqnarray*}
which proves that $(\fg,J,\varphi)$ is again $T^*_0\fk$.

Let us now suppose that $w_0\ne 0$. Choose an element $u_1\in\mbox{Ker}(F)$ such that 
$\varphi_0(u_1,w_0)=1$ and $\varphi_0(Ju_1,w_0)=0$. Take $\beta_1,\beta_2\in\R$ such that $t_0=\beta_1u_1+\beta_2Ju_1$ and consider $u_2=\beta_2Ju_1+w_0$. One inmediately gets $\varphi(u_2,u_2)=\varphi(u_2,Ju_1)=0$ and $\varphi(u_1,u_2)=1$. Further, $s_0=\beta_1u_1+u_2$. From Lemma \ref{algeb2} and Proposition \ref{pR2} we now get
\begin{eqnarray*}
&  [v,Jv]=u_2+\beta_1u_1,\,\,\, [v,u_1]=-Jz,\,\,\, [v,u_2]=aJu_1-\beta_1Jz,\,\,\, [v,Ju_2]=-au_1&\\
&  [Jv,u_1]=z,\,\,\, [Jv,u_2]=bJu_1+\beta_1z,\,\,\, [Jv,Ju_2]=-bu_1,\,\,\, [u_2,Ju_2]=az+bJz.&
\end{eqnarray*}
If we choose in this case
\begin{eqnarray*}
& & x_1={(a^2+b^2)^{-3/5}}(av+bJv),\quad x_3=(a^2+b^2)^{-1/5}u_2,\\
& & x_1^*={(a^2+b^2)^{-2/5}}(az+bJz)\,\quad x_3^*=(a^2+b^2)^{1/5}u_1\\
& &x_2=Jx_1,\quad x_4=Jx_3,\quad x_2^*=Jx_1^*,\quad x_4^*=Jx_3^*,
\end{eqnarray*} 
and we take $\alpha=(a^2+b^2)^{-2/5}\beta_1$, we get that the brackets in $\fg$ are given by
\begin{eqnarray*}
& & [x_1,x_2]=x_3+\alpha x_3^*,\quad [x_1,x_3]=-\alpha x^*_2+x_4^*,\\
& & [x_1,x_4]=-x^*_3,\quad [x_2,x_3]=\alpha x^*_1,\\
& &  [x_3,x_4]=x^*_1,\quad
[x_3^*,x_1]=x_2^*,\quad [x_3^*,x_2]=-x^*_1.\end{eqnarray*} 
Besides, one easily shows that we also have
$$ \varphi(x_i,x_j^*)=\delta_{ij},\quad \varphi(x_i,x_j)=\varphi(x_i^*,x_j^*)=0,\quad 1\le i,j\le 4,$$
and, hence, we get the pseudo-Hermitian structure of $T^*_\theta\fk$ for $\theta=\theta_3+\alpha \theta_1$.

Now, consider the  vectors
\begin{eqnarray*}
& & y_1=x_1+\alpha x_3, \quad y_2=x_2+\alpha x_4,\quad y_3=x_3,\quad y_4=x_4,\\
& & y_1^*=x_1^*, y_2^*=x_2^*,\quad y_3^*=x_3^*-\alpha x_1^*,\quad 
y_4^*=x_4^*-\alpha x_2^*.
\end{eqnarray*}

Notice that $Jy_1=y_2$, $Jy_3=y_4$, $Jy_1^*=y^*_2$, $Jy_3^*=y_4^*$ and a direct computation proves that
$$ \varphi(y_i,y_j^*)=\delta_{ij},\quad \varphi(y_i,y_j)=\varphi(y_i^*,y_j^*)=0,\quad 1\le i,j\le 4.$$
Moreover, $y_1^*,y_2^*,y_4^*$ are clearly central and we have
\begin{eqnarray*}
& & [y_1,y_2]=[x_1,x_2]+\alpha [x_3,x_2]+\alpha [x_1,x_4]+\alpha^2[x_3,x_4]=x_3+\alpha x_3^*-\alpha^2 x^*_1-\alpha x^*_3+\alpha^2 x^*_1=y_3\\
& & [y_1,y_3]=[x_1,x_3]=-\alpha x^*_2+x_4^*=y_4^*\\
& & [y_1,y_4]=[x_1,x_4]+\alpha[x_3,x_4]=-x^*_3+\alpha x^*_1=-y^*_3\\
& & [y_2,y_3]=[x_2,x_3]+\alpha [x_4,x_3]=\alpha x^*_1-\alpha x^*_1=0\\
& & [y_3,y_4]=[x_3,x_4]=x^*_1=y_1^*\\
& & [y_3^*,y_1]=[x_3^*,x_1]+\alpha[x_3^*,x_3]=x_2^*=y_2^*\\
& & [y_3^*,y_2]=[x_3^*,x_2]+\alpha[x_3^*,x_4]=-x_1^*=-y_1^*,
\end{eqnarray*}
and the remaining brackets vanish. So, we actually have the Lie algebra $T^*_{\theta_3}\fk$.\qed 

\bigskip

We can now summarize all the results above. We will denote by ${\mathcal L}^{4,2}$ the Lie algebra with complex structure ${\mathcal L}$ of Corollary \ref{LH} endowed with the Lorentz-Hermitian metric and by ${\mathcal L}^{2,4}$ the same Lie algebra with the opposite metric. As before, the metrics and complex structures on $T^*_\theta \fk$ are the ones given en Proposition \ref{kodaira}.

Since the direct sum of pseudo-Hermitian quadratic Lie algebras is clearly a pseudo-Hermitian quadratic Lie algebra, we will only give the classification of the {\it indeocmposable} ones, this is to say, those pseudo-Hermitian quadratic Lie algebras that do not split as a sum of two  pseudo-Hermitian quadratic ideals.

\begin{theor} \label{c8} For $n\le 8$ a $n$-dimensional indecomposable pseudo-Hermitian quadratic Lie algebra is equivalent to one and only one of the following pseudo-Hermitian  Lie algebras: $\R^{0,2}$, $\R^{2,0}$, ${\mathcal L}^{4,2}$, ${\mathcal L}^{2,4}$, $T^*_0\fk$, $T^*_{\theta_3}\fk$.
\end{theor} 
\proof The result follows at once from the propositions above. It only remains to prove that those Lie algebras are inequivalent. Using the dimension of the Lie algebra and the signature of the metric, it suffices to see that the two T$^*$-extensions are inequivalent. But the center in $T^*_0\fk$ is 5-dimensional whereas that of $T^*_{\theta_3}\fk$ is 3-dimensional, so that they are not isomorphic (actually,  one has that $T^*_0\fk$ is 2-step nilpotent while the other one is 3-step nilpotent).\qed

\bigskip

It is clear that  abelian psuedo-Hermitian quadratic Lie algebras are  classified by the dimension of the Lie algebra and the signature of the metric. The following table gives a classification scheme for non-abelian pseudo-Hermitian quadratic Lie algebras $(\fg, J,\varphi)$ with $\dim(\fg)\le 8$ in terms of the dimensions of $\fg$ and $[\fg,\fg]$ and the signatures of $\varphi$ and its restriction $\varphi_{[\fg,\fg]}$ to $[\fg,\fg]\times [\fg,\fg]$. Although it is not necessary, we include a column with the nilpotency index of $\fg$.

\smallskip

\begin{center}
\begin{tabular}{|c|c|c|c|c|c|}
\hline $\dim (\fg)$&$\dim([\fg,\fg])$&$\mbox{sig}(\varphi)$&$\mbox{sig}(\varphi_{[\fg,\fg]})$& Nilp. Index&$\begin{array}{c} \,\\ \,\end{array}\fg\begin{array}{c} \,\\ \,\end{array}$\\ \hline
6& 3&(2,4)&(0,1)&3&$\begin{array}{c} \,\\ \,\end{array}{\mathcal L}^{2,4}\begin{array}{c} \,\\ \,\end{array}$\\
\hline
6& 3&(4,2)&(1.0)&3&$\begin{array}{c} \,\\ \,\end{array}{\mathcal L}^{4,2}\begin{array}{c} \,\\ \,\end{array}$\\
\hline
8& 3&(2,6)&(0,1)&3&$\begin{array}{c} \,\\ \,\end{array}{\mathcal L}^{2,4}\oplus\R^{0,2}\begin{array}{c} \,\\ \,\end{array}$\\
\hline
8& 3&(4,4)&(0,1)&3&$\begin{array}{c} \,\\ \,\end{array}{\mathcal L}^{2,4}\oplus\R^{2,0}=T^*_{\theta_1}\fk\begin{array}{c} \,\\ \,\end{array}$\\
\hline
8& 3&(4,4)&(0,0)&2&$\begin{array}{c} \,\\ \,\end{array}T^*_{0}\fk\begin{array}{c} \,\\ \,\end{array}$\\
\hline
8& 3&(4,4)&(1,0)&3&$\begin{array}{c} \,\\ \,\end{array}{\mathcal L}^{4,2}\oplus\R^{0,2}=T^*_{-\theta_1}\fk\begin{array}{c} \,\\ \,\end{array}$\\
\hline
8& 3&(6,2)&(1,0)&3&$\begin{array}{c} \,\\ \,\end{array}{\mathcal L}^{4,2}\oplus\R^{2,0}\begin{array}{c} \,\\ \,\end{array}$\\
\hline
8& 5&(4,4)&(1,1)&3&$\begin{array}{c} \,\\ \,\end{array}T^*_{\theta_3}\fk\begin{array}{c} \,\\ \,\end{array}$\\
 \hline
\end{tabular}

\medskip

Table I. {\it Non-abelian  pseudo-Hermitian quadratic nilpotent Lie algebras up to dimension 8}

\end{center}

\bigskip

\noindent {\it Authors' addresses:}

\medskip

\noindent Mustapha Bachaou, Department of Mathematics, Faculty of Sciences, University  Abdelmalek Essa\^adi, B.P. 2121,Tetouan, Morocco. 
\\
mustapha.bachaou@etu.uae.ac.ma

\medskip

\noindent Ignacio Bajo, Departamento de Matem\'atica Aplicada II, E.E. Telecomunicaci\'on, Universidade de Vigo, 36280 Vigo, Spain.
\\
ibajo@dma.uvigo.es

\medskip

\noindent Mohamed Louzari, Department of Mathematics, Faculty of Sciences, University  Abdelmalek Essa\^adi, B.P. 2121,Tetouan, Morocco. 
\\
m.louzari@uae.ac.ma


\begin{thebibliography}{99}

\bibitem{Andrada1} A. Andrada, M.L. Barberis, I.G. Dotti, Abelian Hermitian geometry, {\it  Diff. Geom. Appl.} 30 (2012), no. 5, 509--519.

\bibitem{Andrada} A. Andrada, M.L. Barberis, G. Ovando, Lie bialgebras of complex type and associated Poisson Lie groups, {\it Jour. Geom. Phys.} 58 (2008), 1310--1328.

\bibitem{bajo1} I. Bajo, Prederivations of Lie algebras and isometries
of bi-invariant Lie groups, {\it Geom. Dedicata} 66 (1997), no. 3,  281--291.

\bibitem{IB} I. Bajo, S. Benayadi, A. Medina, Symplectic structures on quadratic Lie algebras, {\it Jour. Algebra} 316 (2007), 174--188.

\bibitem{BD} M.L. Barberis, I. Dotti, Abelian complex structures on solvable Lie algebras, {\it J. Lie Theory} 14 (2004), 25--34. 

\bibitem{Borde} M. Bordemann, Nondegenerate invariant bilinear forms on nonassociative algebras, {\it Acta Math. Univ. Com.}, LXVI, 2 (1997),
151--201.

\bibitem{FS} G. Favre, L.J. Santharoubane, Symmetric, invariant, Non-degenerate
Bilinear Form on a Lie Algebra, \textit{Jour. Algebra} 105, 451-464 (1987).

\bibitem{hk} K.H. Hofmann, V.S. Keith,
Invariant quadratic forms on finite-dimensional Lie algebras,
{\it Bull. Austral. Math. Soc.} 33 (1986), no. 1, 21-–36. 

\bibitem{IK} I. Kath, Nilpotent metric Lie algebras of small dimension, \textit{J. Lie Theory}
 17 (2007) 41-–61.

\bibitem{Koda} K. Kodaira, On the structure of compact complex analytic surfaces I, \textit{Amer. J. Math.} 86 (1964), 751--798.


\bibitem{MR} A. Medina, Ph. Revoy, Algèbres de Lie et produit scalaire invariant, \textit{Ann. Scient. l'É.N.S.} 4$^e$ série, t. 18, 3 (1985),  553--561.


\bibitem{Million} D.V. Millionshchikov, Complex structures on nilpotent Lie algebras and descending central series, \textit{Rend. Sem. Mat. Univ. Pol. Torino} 74, 1 (2016), 163--182.

\bibitem{salamon} S.M. Salamon, Complex structures on nilpotent Lie algebras, \textit{Jour. Pure  App. Algebra} 157 (2001), 311--333.

\bibitem{Samelson} H. Samelson, A class of compact-analytic manifolds, \textit{Portug. Math.} 12 (1953),  129--132.
 
\bibitem{Thurs} W. Thurston, Some simple examples of symplectic manifolds, \textit{Proc. Amer. Math. Soc.} 55 no. 2 (1976), 467--468.

\end{thebibliography}
\end{document}